\documentclass[12pt,oneside,leqno]{amsproc}
\usepackage[utf8]{inputenc}
\usepackage[russian]{babel}
\usepackage{amssymb,amsmath,amsthm}
\renewcommand{\subsection}{\refstepcounter{subsection}%
\par\bigskip\noindent\textbf{\upshape\thesubsection. }}
\renewcommand{\subsubsection}{\refstepcounter{subsubsection}%
\par\medskip\noindent\textbf{\upshape\thesubsubsection.  }}
\renewcommand{\paragraph}{\refstepcounter{paragraph}%
\par\smallskip\noindent\textbf{\upshape\theparagraph. }}

\numberwithin{equation}{subsection}

\renewcommand{\thesubsection}{\arabic{subsection}}
\renewcommand{\thesubsubsection}{\arabic{subsection}.\arabic{subsubsection}}

\newcommand{\Wo}{{\raisebox{0.2ex}{\(\stackrel{\circ}{W}\)}}{}}

\makeatletter
\renewcommand{\@makefnmark}{\hbox{\small\mathsurround=0cm%
${}\hspace{0.04cm}{}^{\@thefnmark)}$}}
\renewcommand{\@makefntext}[1]{\parindent=1em\noindent\hbox to 1.8em{%
\hss${}^{\@thefnmark)}$}\,#1}
\makeatother

\title{О задаче Неймана для уравнения Штурма--Лиувилля с самоподобным весом
канторовского типа}
\author{А.~А.~Владимиров, И.~А.~Шейпак\footnote{Работа поддержана РФФИ,
грант \No~10-01-00423.}}
\begin{document}
\renewcommand{\proofname}{{\upshape Д\,о\,к\,а\,з\,а\,т\,е\,л\,ь\,с\,т\,в\,о.}}
\begin{abstract}
Рассматриваются вторая и третья граничные задачи для дифференциального уравнения
\[
    -y''-\lambda\rho y=0,
\]
где вес \(\rho\in W_2^{-1}[0,1]\) представляет собой обобщённую производную
самоподобной функции канторовского типа. На основе изучения осцилляционных свойств
собственных функций существенно уточняются характеристики известных спектральных
асимптотик таких задач.
\end{abstract}
\begin{flushleft}
\normalsize УДК~517.984
\end{flushleft}
\maketitle
\markboth{}{}

\section{Введение}\label{par:1}
\subsection\label{pt:1}
В настоящей статье продолжается начатое работами \cite{SV} и \cite{VSh1}
изучение спектральных асимптотик граничной задачи
\begin{gather}\label{eq:1.1}
    -y''-\lambda\rho y=0,\\ \label{eq:1.2}
    y(0)=y(1)=0,
\end{gather}
где вес \(\rho\in W_2^{-1}[0,1]\) представляет собой обобщённую производную
некоторой самоподобной функции \(P\in L_2[0,1]\). Основные результаты указанных
работ состоят в следующем. Установлено, что в случае строгой положительности
\emph{спектрального порядка} \(D\in [0,1)\) функции \(P\) считающая функция
\(N:(0,+\infty)\to\mathbb N\) собственных значений задачи \eqref{eq:1.1},
\eqref{eq:1.2} имеет при \(\lambda\to+\infty\) асимптотику
\begin{equation}\label{eq:mes_asymp}
	N(\lambda)=\lambda^D\cdot\bigl(s(\ln\lambda)+o(1)\bigr),
\end{equation}
где \(s\) "--- некоторая зависящая от выбора веса \(\rho\) непрерывная
периодическая функция. В случае \emph{неарифметичности} типа самоподобия функции
\(P\) функция \(s\) заведомо является постоянной. В случае же \emph{арифметичности}
этого типа функция \(s\) обладает периодом, равным \emph{шагу самоподобия}
функции \(P\).

В работе \cite{Naz} была сформулирована гипотеза о том, что функция \(s\) является
непостоянной для произвольного неравномерного веса с арифметически самоподобной
первообразной. В работе \cite{VSh1} на основе результатов численных экспериментов
было показано, что функция \(s\) действительно не может являться постоянной в том
простейшем случае, когда обобщённая первообразная \(P\) веса \(\rho\) представляет
собой канторову лестницу. Одной из основных целей настоящей статьи является
проводимое традиционными теоретическими методами установление справедливости
указанной гипотезы (и даже некоторых более сильных утверждений) для произвольной
функции \(P\), имеющей \emph{канторовский тип} самоподобия. Точные определения
соответствующих понятий будут даны нами далее.

Следует отметить, что в случае \(D=0\) асимптотика собственных значений задачи
\eqref{eq:1.1}, \eqref{eq:1.2} имеет вид, коренным образом отличный от
\eqref{eq:mes_asymp} (см. \cite{VSh2010}). Однако для целей настоящей статьи
это обстоятельство является безразличным.

\subsection\label{pt:1.2}
В отличие от работ \cite{SV} и \cite{VSh1}, в настоящей статье мы будем опираться
не на теорию восстановления, а на осцилляционную теорию для задач Штурма--Лиувилля
с сингулярными коэффициентами, развитую в недавних работах \cite{BASh2009}
и \cite{V2}\footnote{Отдельного упоминания заслуживает также работа \cite{BNT},
в которой были получены близкие к используемым далее результаты о нулях собственных
функций для фрактальных операторов Штурма--Лиувилля.}. При этом основной упор будет
делаться нами на изучение осцилляционных свойств собственных функций не первой,
а второй и третьей граничных задач для уравнения \ref{pt:1}\,\eqref{eq:1.1}.
Именно для таких задач будет установлено наличие эффекта \emph{спектральной
периодичности}\footnote{Оно было анонсировано в тезисах \cite{VSh2004}.},
играющего ключевую роль при последующем уточнении характеристик асимптотики
\ref{pt:1}\,\eqref{eq:mes_asymp}. Используемая далее терминология соответствует
результатам заметки \cite{Vl}, согласно которым поставленной в терминах
\emph{квазипроизводных} \cite[\S~15]{Na} формальной граничной задаче
\begin{gather}\label{eq:rp1}
	-y''-\lambda\rho y=0,\\ \label{eq:rp2}
	y^{[1]}(0)-\gamma_0 y(0)=y^{[1]}(1)+\gamma_1 y(1)=0
\end{gather}
сопоставляется линейный операторный пучок \(T_{\rho}:\mathbb C\to\mathcal B(
W_2^1[0,1],W_2^{-1}[0,1])\) со свойством
\[
	(\forall\lambda\in\mathbb C)\,(\forall y\in W_2^1[0,1])\qquad
	\langle T_{\rho}(\lambda)y,y\rangle=\int\limits_0^1 |y'|^2\,dx-
	\lambda\langle\rho,|y|^2\rangle+\gamma_0\,|y(0)|^2+\gamma_1\,|y(1)|^2.
\]

\subsection
Структура статьи имеет следующий вид. В \ref{par:sfkt} приводятся необходимые
для дальнейшего св\'{е}дения о самоподобных функциях канторовского типа и связанных
с ними сингулярных граничных задачах. В \ref{par:2} устанавливается факт
спектральной периодичности для некоторых из таких граничных задач. В \ref{par:5},
носящем вспомогательный характер, устанавливается применяемый в дальнейшем
критерий сингулярности монотонной функции. В \ref{par:3} проводится основанное
на использовании спектральной периодичности уточнение характеристик коэффициента
\(s\) из асимптотики \ref{pt:1}\,\eqref{eq:mes_asymp}. Наконец, в \ref{par:4}
приводятся иллюстрирующие явление спектральной периодичности расчётные данные.

Поясняющие изложение ссылки приводятся в квадратных скобках. В общем случае они
содержат указание на номер параграфа и пункта. В случае, когда ссылка даётся внутри
одного параграфа, указание его номера опускается.

Нумерация формул подчиняется нумерации пунктов. В случае, когда ссылка на формулу
даётся внутри одного пункта, указание номера этого пункта опускается.


\section{Самоподобные функции канторовского типа}\label{par:sfkt}
\subsection\label{pt:sfkt}
Пусть фиксировано некоторое натуральное число \(\varkappa>1\), а также
два вещественных числа \(a\in(0,1/\varkappa)\) и \(b\rightleftharpoons
(1-\varkappa a)/(\varkappa-1)\). Свяжем с этими числами набор точек \(\alpha_{2k}
\rightleftharpoons k(a+b)\) и \(\alpha_{2k+1}\rightleftharpoons\alpha_{2k}+a\),
где \(k\in\{0,\ldots,\varkappa-1\}\). Функцию \(f\in C[0,1]\), удовлетворяющую
равенствам \(f(0)=0\) и \(f(1)=1\), мы будем называть \emph{самоподобной функцией
канторовского типа} с параметрами \(\varkappa\), \(a\) и \(b\), если при любом
выборе индекса \(k\in\{1,\ldots,\varkappa-1\}\) она постоянна на интервале
\((\alpha_{2k-1},\alpha_{2k})\), и при любом выборе индекса \(k\in\{0,\ldots,
\varkappa-1\}\) функция \(f_k\in C[0,1]\) вида
\[
	f_k(x)\rightleftharpoons\varkappa f(\alpha_{2k}+ax)
\]
совпадает с функцией \(f\) с точностью до аддитивной постоянной.

В частности, классическая канторова лестница представляет собой самоподобную
функцию канторовского типа с параметрами \(\varkappa=2\), \(a=b=1/3\).

На основе известных результатов \cite{VSh1}, \cite{Sh} легко устанавливается,
что каждый допустимый набор параметров \(\varkappa\), \(a\) и \(b\) определяет
единственную самоподобную функцию. Более того, эта функция является арифметически
самоподобной с шагом \(\nu=\ln(\varkappa/a)\), а её спектральный порядок
определяется равенством \(D=\nu^{-1}\ln\varkappa\in (0,1/2)\).

Легко проверить, что для любой самоподобной функции \(f\) канторовского типа
функция \(1-f(1-\cdot)\) представляет собой самоподобную функцию канторовского
типа с теми же параметрами \(\varkappa\), \(a\) и \(b\), что и исходная. Тем самым,
справедливо следующее утверждение:

\subsubsection\label{invert}
{\itshape Пусть \(f\in C[0,1]\) "--- самоподобная функция канторовского типа.
Тогда выполняется тождество
\[
	(\forall x\in [0,1])\qquad f(x)=1-f(1-x).
\]
}

\subsection\label{pt:2.2}
В отличие от случая задачи Дирихле, где вес \(\rho\in\Wo_2^{-1}[0,1]\) полностью
задаётся своей квадратично суммируемой обобщённой первообразной \(P\), в случае
задачи Неймана для определения веса \(\rho\in W_2^{-1}[0,1]\) формулой
\[
	(\forall y\in W_2^1[0,1])\qquad \langle\rho,y\rangle=
	-\int\limits_0^1 P\overline{y'}\,dx+\left.P\overline{y}\right.|_0^1
\]
дополнительно требуется указание "`граничных значений"' \(P(0)\) и \(P(1)\).
Если обобщённая первообразная \(P\) является непрерывной, то на роль таких
значений естественным образом выбираются "`обычные"' граничные значения функции
\(P\). Соответственно, квадратичная форма операторов \(T_{\rho}(\lambda)\)
из пучка, отвечающего граничной задаче \ref{par:1}.\ref{pt:1.2}\,\eqref{eq:rp1},
\ref{par:1}.\ref{pt:1.2}\,\eqref{eq:rp2}, где вес \(\rho\) представляет собой
обобщённую производную самоподобной функции \(P\) канторовского типа, принимает вид
\begin{equation}\label{eq:rp3}
	\langle T_{\rho}(\lambda)y,y\rangle=\int\limits_0^1 \bigl\{|y'|^2+
		\lambda P(|y|^2)'\bigr\}\,dx+\gamma_0\,|y(0)|^2+
		(\gamma_1-\lambda)\,|y(1)|^2.
\end{equation}

Известно \cite{SaSh}, \cite{Vl}, что любая функция \(y\), принадлежащая ядру
определённого тождеством \eqref{eq:rp3} оператора \(T_{\rho}(\lambda)\), обладает
абсолютно непрерывной удлинённой производной \(y^{<1>}\rightleftharpoons y'+
\lambda Py\), удовлетворяющей соотношениям
\begin{gather*}
	(y^{<1>})'=-\lambda^2P^2y+\lambda Py^{<1>},\\
	y^{<1>}(0)-\gamma_0 y(0)=y^{<1>}(1)+(\gamma_1-\lambda) y(1)=0.
\end{gather*}
Ввиду непрерывности функции \(P\), отсюда немедленно вытекает непрерывная
дифференцируемость функции \(y\) и связанная с этим возможность прямого понимания
граничных условий \ref{par:1}.\ref{pt:1.2}\,\eqref{eq:rp2}. Отсюда же может быть
получен \cite[Утверждение~11]{V2} следующий важный для дальнейшего факт:

\subsubsection\label{oscil}
{\itshape Пусть \(\{\lambda_n\}_{n=0}^{\infty}\) "--- последовательность
занумерованных в порядке возрастания собственных значений граничной задачи
\ref{par:1}.\ref{pt:1.2}\,\eqref{eq:rp1}, \ref{par:1}.\ref{pt:1.2}\,%
\eqref{eq:rp2}. Тогда независимо от выбора индекса \(n\in\mathbb N\) собственное
значение \(\lambda_n\) является простым, причём любая отвечающая ему собственная
функция не обращается в нуль на границе отрезка \([0,1]\) и имеет внутри этого
отрезка в точности \(n\) различных нулей.
}

\subsection
Легко видеть [\ref{pt:2.2}\,\eqref{eq:rp3}], что замена коэффициентов \(\gamma_0\)
и \(\gamma_1\) в граничных условиях \ref{par:1}.\ref{pt:1.2}\,\eqref{eq:rp2}
приводит к возмущению каждого из операторов \(T_{\rho}(\lambda)\) некоторым
оператором ранга \(2\). Из общей вариационной теории пучков самосопряжённых
операторов (см., например, \cite{LSY}) потому следует, что считающие функции
собственных значений различных третьих граничных задач, отвечающих одному
и тому же уравнению \ref{par:1}.\ref{pt:1.2}\,\eqref{eq:rp1}, не могут различаться
более, чем на \(2\). Аналогичная связь имеется между считающими функциями третьей
граничной задачи и задачи Дирихле "--- квадратичная форма которой получается
из квадратичной формы третьей граничной задачи сужением на подпространство
коразмерности \(2\). Тем самым, асимптотическая формула \ref{par:1}.\ref{pt:1}\,%
\eqref{eq:mes_asymp} справедлива для спектров всех рассматриваемых далее
граничных задач.


\section{Спектральная периодичность}\label{par:2}
\subsection
Имеют место следующие три факта:

\nopagebreak
\subsubsection\label{3.2.2}
{\itshape Пусть \(\{\lambda_n\}_{n=0}^{\infty}\) "--- последовательность
занумерованных в порядке возрастания собственных значений отвечающей уравнению
\ref{par:1}.\ref{pt:1.2}\,\eqref{eq:rp1} граничной задачи
\begin{equation}\label{eq:2.0}
	y^{[1]}(0)=y^{[1]}(1)=0.
\end{equation}
Тогда независимо от выбора индекса \(n\in\mathbb N\) выполняется равенство
\begin{equation}\label{eq:2.101}
	\lambda_{\varkappa n}=(\varkappa/a)\,\lambda_n.
\end{equation}
}

\begin{proof}
Зафиксируем отвечающую собственному значению \(\lambda_n\) собственную функцию
\(y_n\), не обращающуюся в нуль на границе отрезка \([0,1]\) и имеющую внутри него
в точности \(n\) различных нулей [\ref{par:sfkt}.\ref{oscil}]. Сопоставим ей
нетривиальную функцию \(z\in C[0,1]\), удовлетворяющую следующим условиям:
\begin{enumerate}
\item При любом выборе индекса \(k\in\{1,\ldots,\varkappa-1\}\) функция \(z\)
постоянна на интервале \((\alpha_{2k-1},\alpha_{2k})\).
\item При любом выборе индекса \(k\in\{0,\ldots,\varkappa-1\}\) функция \(z_k\) вида
\[
	z_k(x)\rightleftharpoons z(\alpha_{2k}+ax)
\]
совпадает с функцией \(y_n\) с точностью до мультипликативной постоянной.
\end{enumerate}
С учётом сделанных в пункте \ref{par:sfkt}.\ref{pt:2.2} замечаний, граничные
условия \eqref{eq:2.0} гарантируют непрерывную дифференцируемость функции \(z\).
При этом факт самоподобия функции \(P\) позволяет легко проверить, что
удлинённая производная \(z^{<1>}\rightleftharpoons z'+(\varkappa/a)\,\lambda_n Pz\)
абсолютно непрерывна и удовлетворяет соотношениям
\begin{gather*}
	(z^{<1>})'=-(\varkappa/a)^2\,\lambda_n^2 P^2 z+
		(\varkappa/a)\,\lambda_n Pz^{<1>},\\
	z^{<1>}(0)=z^{<1>}(1)-(\varkappa/a)\,\lambda_n z(1)=0.
\end{gather*}
Тем самым \cite{SaSh}, \cite{Vl}, функция \(z\) представляет собой собственную
функцию граничной задачи \ref{par:1}.\ref{pt:1.2}\,\eqref{eq:rp1}, \eqref{eq:2.0},
отвечающую собственному значению \((\varkappa/a)\,\lambda_n\). При этом она имеет
на интервале \((0,1)\) в точности \(\varkappa n\) различных нулей, что и означает
[\ref{par:sfkt}.\ref{oscil}] выполнение равенства \eqref{eq:2.101}.
\end{proof}

\subsubsection\label{3.2.3}
{\itshape Пусть \(\{\lambda_n\}_{n=0}^{\infty}\) "--- последовательность
занумерованных в порядке возрастания собственных значений отвечающей уравнению
\ref{par:1}.\ref{pt:1.2}\,\eqref{eq:rp1} граничной задачи
\[
	by^{[1]}(0)-2y(0)=by^{[1]}(1)+2y(1)=0,
\]
а \(\{\mu_n\}_{n=0}^{\infty}\) "--- аналогичная последовательность для отвечающей
тому же уравнению граничной задачи
\[
	by^{[1]}(0)-2ay(0)=by^{[1]}(1)+2ay(1)=0.
\]
Тогда независимо от выбора индекса \(n\in\mathbb N\) выполняется равенство
\[
	\lambda_{\varkappa (n+1)-1}=(\varkappa/a)\,\mu_n.
\]
}

Доказательство проводится аналогичным доказательству утверждения \ref{3.2.2}
способом. Основное отличие заключается в том, что используемые при "`сшивке"'
копий исходной собственной функции горизонтальные отрезки заменяются на наклонные
с нулями в центрах соответствующих интервалов \((\alpha_{2k-1},\alpha_{2k})\).

\subsubsection\label{3.2.100}
{\itshape Пусть \(\{\lambda_n\}_{n=0}^{\infty}\) "--- последовательность
занумерованных в порядке возрастания собственных значений граничной задачи
\ref{par:1}.\ref{pt:1.2}\,\eqref{eq:rp1}, \eqref{eq:2.0}, а \(\{\mu_n\}_{n=0}^{%
\infty}\) "--- аналогичная последовательность для отвечающей тому же уравнению
граничной задачи
\begin{equation}\label{eq:2.102}
	y^{[1]}(0)=by^{[1]}(1)+2ay(1)=0.
\end{equation}
Тогда в случае чётности параметра \(\varkappa\) независимо от выбора индекса
\(n\in\mathbb N\) выполняется равенство
\[
	\lambda_{\varkappa (n+1/2)}=(\varkappa/a)\,\mu_n.
\]
}

Доказательство также проводится аналогичным доказательству утверждения \ref{3.2.2}
способом. При этом учитывается то обстоятельство [\ref{par:sfkt}.\ref{invert},
\ref{par:sfkt}.\ref{pt:2.2}\,\eqref{eq:rp3}], что для любой собственной функции
\(y_n\) граничной задачи \ref{par:1}.\ref{pt:1.2}\,\eqref{eq:rp1}, \eqref{eq:2.102}
функция \(y_n(1-\cdot)\) представляет собой собственную функцию отвечающей тому же
уравнению граничной задачи
\[
	by^{[1]}(0)-2ay(0)=y^{[1]}(1)=0.
\]
Собственная функция \(z\) граничной задачи \ref{par:1}.\ref{pt:1.2}\,%
\eqref{eq:rp1}, \eqref{eq:2.0} составляется при этом из чередующихся копий функций
\(y_n\) и \(y_n(1-\cdot)\), "`сшиваемых"' наклонными и горизонтальными отрезками
попеременно.


\section{Критерий сингулярности}\label{par:5}
\subsection
На протяжении настоящего параграфа мы отступаем от проводимого в остальной части
статьи закрепления символов \(a\) и \(b\) за параметрами самоподобия.

\subsubsection\label{vspom}
{\itshape Пусть \(f\in L_2[a,b]\) "--- неубывающая функция, а \(A\) и \(B\) "---
два различных вещественных числа, почти всюду на отрезке \([a,b]\) удовлетворяющие
неравенствам
\[
	A\leqslant f(x)\leqslant B.
\]
Тогда при любых \(n\in\mathbb N\) и \(\varepsilon>0\) найдётся ступенчатая
функция \(f_n\) с не более чем \(2^n\) точками разрыва на интервале \((a,b)\),
тождественно равная \(A\) на некоторой правой окрестности точки \(a\), тождественно
равная \(B\) на некоторой левой окрестности точки \(b\), и удовлетворяющая
неравенству
\[
	\|f-f_n\|_{L_2[a,b]}<2^{-n}\cdot(1+\varepsilon)\cdot(B-A)\cdot\sqrt{b-a}.
\]
}

\begin{proof}
Без ограничения общности рассмотрения можно считать функцию \(f\) непрерывной.
Это и будет далее предполагаться.

Воспользуемся методом арифметической индукции. При \(n=0\) на роль искомой функции
\(f_0\) может быть взята функция, принимающая значение \(A\) слева от точки
\((a+b)/2\), и значение \(B\) "--- справа. Далее, зафиксируем произвольное число
\(\delta>0\), удовлетворяющее неравенству
\begin{equation}\label{eq:100}
	(1+\delta)^2<1+\varepsilon.
\end{equation}
Заведомо имеет место один из следующих трёх случаев:
\begin{enumerate}
\item Найдётся точка \(\zeta\in (a,b)\), удовлетворяющая неравенству
\[
	\left|f(\zeta)-\dfrac{B+A}{2}\right|<\delta\cdot \dfrac{B-A}{2}.
\]
\item Выполняется неравенство
\[
	f(a)>\dfrac{B+A}{2}-\delta\cdot\dfrac{B-A}{2}.
\]
\item Выполняется неравенство
\[
	f(b)<\dfrac{B+A}{2}+\delta\cdot\dfrac{B-A}{2}.
\]
\end{enumerate}

В первом случае, согласно индуктивному предположению, найдётся ступенчатая
функция \(f_{n+1}\), имеющая не более \(2^n\) точек разрыва на каждом из
интервалов \((a,\zeta)\) и \((\zeta,b)\), тождественно равная \(A\) на некоторой
правой окрестности точки \(a\), тождественно равная \(B\) на некоторой левой
окрестности точки \(b\), принимающая значение \(f(\zeta)\) в некоторой окрестности
точки \(\zeta\), а также удовлетворяющая неравенствам
\begin{gather*}
	\|f-f_{n+1}\|_{L_2[a,\zeta]}<2^{-n}\cdot(1+\delta)\cdot (f(\zeta)-A)\cdot
		\sqrt{\zeta-a},\\
	\|f-f_{n+1}\|_{L_2[\zeta,b]}<2^{-n}\cdot(1+\delta)\cdot (B-f(\zeta))\cdot
		\sqrt{b-\zeta}.
\end{gather*}
Ввиду \eqref{eq:100}, последнее влечёт справедливость искомого неравенства
\begin{equation}\label{eq:101}
	\|f-f_{n+1}\|_{L_2[a,b]}<2^{-n-1}\cdot(1+\varepsilon)\cdot(B-A)\cdot
		\sqrt{b-a}.
\end{equation}

Во втором случае, согласно индуктивному предположению, найдётся ступенчатая
функция \(f_n\), имеющая не более \(2^n\) точек разрыва на интервале \((a,b)\),
тождественно равная \(f(a)\) на некоторой правой окрестности точки \(a\),
тождественно равная \(B\) на некоторой левой окрестности точки \(b\), а также
удовлетворяющая неравенству
\begin{gather*}
	\|f-f_n\|_{L_2[a,b]}<2^{-n}\cdot(1+\delta)\cdot (B-f(a))\cdot\sqrt{b-a},
	\intertext{а потому и неравенству}
	\|f-f_n\|_{L_2[a,b]}<2^{-n-1}\cdot(1+\varepsilon)\cdot (B-A)\cdot
		\sqrt{b-a}.
\end{gather*}
Строгий характер последнего неравенства позволяет теперь, изменяя функцию \(f_n\)
на достаточно малой окрестности точки \(a\), построить функцию \(f_{n+1}\),
имеющую не более \(2^n+1\) точек разрыва на интервале \((a,b)\), тождественно
равную \(A\) на некоторой правой окрестности точки \(a\), тождественно равную
\(B\) на некоторой левой окрестности точки \(b\), и удовлетворяющую неравенству
\eqref{eq:101}.

Третий случай разбирается аналогично второму.
\end{proof}

\subsubsection
{\itshape Пусть \(f\in L_2[0,1]\) "--- ограниченная неубывающая чисто сингулярная
функция. Тогда найдётся последовательность \(\{f_n\}_{n=0}^{\infty}\) неубывающих
ступенчатых функций, удовлетворяющая при \(n\to\infty\) асимптотическому
соотношению
\[
	(\#\mathfrak A_n+2)\cdot\|f-f_n\|_{L_2[0,1]}=o(1),
\]
где через \(\{\mathfrak A_n\}_{n=0}^{\infty}\) обозначена последовательность
множеств точек разрыва функций \(f_n\).
}

\begin{proof}
Зафиксируем произвольное вещественное число \(\varepsilon>0\), и представим
\cite[Теорема~13.3]{Saks:1949} сингулярную функцию \(f\) в виде суммы двух
ограниченных неубывающих функций \(g\) и \(h\) следующего вида:
\begin{enumerate}
\item Найдётся вещественное число \(C\), почти всюду на отрезке \([0,1]\)
удовлетворяющее соотношениям
\[
	C\leqslant g(x)\leqslant C+\varepsilon.
\]
\item Найдётся система \(\{\Gamma_k\}_{k=1}^m\) попарно непересекающихся
интервалов, содержащая все точки роста функции \(h\), и имеющая не превосходящую
\(\varepsilon^2\) суммарную длину.
\end{enumerate}
Обозначим через \(a_k\) длины интервалов \(\Gamma_k\), а через \(d_k\) "---
приращения функции \(h\) на этих интервалах. Согласно утверждению \ref{vspom},
при любом \(n\in\mathbb N\) найдётся ступенчатая функция \(h_n\) со следующими
свойствами:
\begin{enumerate}
\item Все точки разрыва функции \(h_n\) содержатся в объединении интервалов
\(\Gamma_k\).
\item На каждом из интервалов \(\Gamma_k\), где \(k\in\{1,\ldots, m\}\),
функция \(h_n\) имеет не более \(2^n d_k+1\) точек разрыва.
\item Выполняются неравенства
\begin{flalign*}
	&& \|h-h_n\|_{L_2[0,1]}^2&\leqslant\sum\limits_{k=1}^m
		(2^{n-1}d_k+1/2)^{-2}\cdot 4a_k\,d_k^2\\
	&& &\leqslant 4^{2-n}\cdot\sum\limits_{k=1}^m a_k\\
	&& &\leqslant 4^{2-n}\varepsilon^2.&&
\end{flalign*}
\end{enumerate}
Аналогичным образом, найдётся имеющая на интервале \((0,1)\) не более \(2^n\)
точек разрыва ступенчатая функция \(g_n\) со свойством
\[
	\|g-g_n\|_{L_2[0,1]}<2^{1-n}\varepsilon.
\]
При этом функция \(f_n\rightleftharpoons h_n+g_n\) имеет не более
\[
	2^n\cdot\left(1+\sum\limits_{k=1}^m d_k\right)+m
\]
точек разрыва на интервале \((0,1)\) и аппроксимирует в смысле среднего
квадратичного функцию \(f\) с точностью не хуже \(2^{3-n}\varepsilon\).
Соответственно, при достаточно больших значениях \(n\) произведение увеличенного
на \(2\) числа точек разрыва функции \(f_n\) и величины \(\|f-f_n\|_{L_2[0,1]}\)
заведомо допускает оценку
\[
	8\varepsilon\cdot(2+B-A),
\]
где \(A\) и \(B\) "--- произвольным образом фиксированные нижняя и верхняя
существенные грани функции \(f\). Ввиду произвольности выбора величины
\(\varepsilon>0\), полученные результаты означают справедливость доказываемого
утверждения.
\end{proof}

\subsubsection\label{asymp}
{\itshape Пусть \(f\in L_2[0,1]\) "--- ограниченная неубывающая функция,
\(\{f_n\}_{n=0}^{\infty}\) "--- последовательность неубывающих ступенчатых
функций, а \(\{\mathfrak A_n\}_{n=0}^{\infty}\) "--- последовательность множеств
точек разрыва функций \(f_n\). Пусть также при \(n\to\infty\) выполняется
асимптотическое соотношение
\[
	(\#\mathfrak A_n+2)\cdot\|f-f_n\|_{L_2[0,1]}=o(1).
\]
Тогда монотонная функция \(f\) является чисто сингулярной.
}

\begin{proof}
Без ограничения общности рассмотрения можно считать, что при любом выборе индекса
\(n\in\mathbb N\) выполняются соотношения
\begin{gather}\notag
	\|f-f_n\|_{L_2[0,1]}\neq 0,\\ \label{eq:3.1}
	2(n+1)^4\cdot(\#\mathfrak A_n+2)\cdot\|f-f_n\|_{L_2[0,1]}<1.
\end{gather}
Это и будет далее предполагаться.

Введём в рассмотрение функции \(\varphi_n:[0,1]\to[0,1]\) вида
\begin{equation}\label{eq:3.2}
	\varphi_n(x)\rightleftharpoons\inf\bigl(1,\,\|f-f_n\|_{L_2[0,1]}^{-1}
	\mathop{\operatorname{dist}}(x,\mathfrak A_n\cup\{0,1\})\bigr).
\end{equation}
Независимо от выбора индекса \(n\in\mathbb N\) функция \(\varphi_n\) является
абсолютно непрерывной, и абсолютная величина её производной \(\varphi'_n\) почти
всюду мажорируется величиной \(\|f-f_n\|_{L_2[0,1]}^{-1}\). Кроме того, при любом
\(n\in\mathbb N\) функция \(\varphi_n\) тождественно равна \(1\) на некотором
открытом множестве, имеющем дополнение не превосходящей \((n+1)^{-4}\) меры
[\eqref{eq:3.1}]. Соответственно, выполняются соотношения
\begin{flalign*}
	&& \int\limits_0^1\varphi_n\,df&=\int\limits_0^1\varphi_n\,df_n+
		\int\limits_0^1\varphi_n\,d(f-f_n)\\
	&& &=\int\limits_0^1\varphi_n\,d(f-f_n)&\text{[\eqref{eq:3.2}]}&\\
	&& &=-\int\limits_0^1(f-f_n)\varphi'_n\,dx&\text{[\eqref{eq:3.2}]}&\\
	&& &\leqslant\|f-f_n\|_{L_2[0,1]}\cdot\|\varphi'_n\|_{L_2[0,1]}\\
	&& &\leqslant (n+1)^{-2},\\
	&& \int\limits_0^1(1-\varphi_n)\,dx&\leqslant (n+1)^{-4},
\end{flalign*}
означающие сходимость числовых рядов
\begin{gather*}
	\sum\limits_{n=0}^{\infty}\int\limits_0^1\varphi_n\,df,\\
	\sum\limits_{n=0}^{\infty}\int\limits_0^1(1-\varphi_n)\,dx.
\end{gather*}
Тогда, согласно теореме~Б.~Леви, множество точек сходимости функционального ряда
\[
	\sum\limits_{n=0}^{\infty}\varphi_n
\]
пренебрежимо по линейной мере Лебега и обладает дополнением, пренебрежимым по мере
с функцией распределения \(f\).
\end{proof}


\section{Уточнение характеристик спектра}\label{par:3}
\subsection
Основной целью настоящего параграфа является установление следующих двух фактов:

\nopagebreak
\subsubsection\label{sing}
{\itshape Коэффициент \(s\) из асимптотики \ref{par:1}.\ref{pt:1}\,%
\eqref{eq:mes_asymp} допускает представление
\[
	(\forall t\in [0,\nu])\qquad s(t)=e^{-Dt}\,\sigma(t),
\]
где \(\sigma\) "--- некоторая чисто сингулярная неубывающая функция.
}

\subsubsection\label{nekonst}
{\itshape Коэффициент \(s\) из асимптотики \ref{par:1}.\ref{pt:1}\,%
\eqref{eq:mes_asymp} не является постоянной функцией.
}

\bigskip
Утверждение \ref{nekonst} представляет собой тривиальное следствие утверждения
\ref{sing}. Действительно, постоянность функции \(s\) равносильна выполнению
на отрезке \([0,\nu]\) тождества \(\sigma(t)\equiv\mathrm{const}\cdot e^{Dt}\),
означающего абсолютную непрерывность функции \(\sigma\). Таким образом,
утверждение \ref{sing} может быть рассмотрено в качестве единственного
центрального результата настоящего параграфа.

\subsection\label{pt:3.1}
При доказательстве утверждения \ref{sing} нами будет использоваться следующий факт:

\nopagebreak
\subsubsection\label{3.2.1}
{\itshape Пусть \(\{\lambda_n\}_{n=0}^{\infty}\) "--- последовательность
занумерованных в порядке возрастания собственных значений граничной задачи
\begin{gather}\label{eq:3.0}
	-y''-\lambda\rho y=0,\\ \label{eq:3.00}
	y^{[1]}(0)=y^{[1]}(1)=0,
\end{gather}
а \(\{\mu_n\}_{n=0}^{\infty}\) "--- аналогичная последовательность для отвечающей
тому же уравнению граничной задачи
\[
	y^{[1]}(0)-\gamma_0y(0)=y^{[1]}(1)+\gamma_1y(1)=0,
\]
где \(\gamma_0,\,\gamma_1\geqslant 0\). Тогда последовательность частичных сумм
числового ряда
\[
	\sum\limits_{n=1}^{\infty} |\ln\mu_n-\ln\lambda_n|
\]
является ограниченной.
}

\begin{proof}
Зафиксируем произвольное значение \(\varepsilon>0\). Заметим, что квадратичная
форма \ref{par:sfkt}.\ref{pt:2.2}\,\eqref{eq:rp3} оператора \(T_{\rho}(
-\varepsilon)\), отвечающего граничной задаче \eqref{eq:3.0}, \eqref{eq:3.00},
является равномерно положительной. Поэтому оператор \(JT_{\rho}(-\varepsilon):
W_2^1[0,1]\to W_2^1[0,1]\), где через \(J:W_2^{-1}[0,1]\to W_2^1[0,1]\) обозначена
изометрия
\[
	(\forall y\in W_2^{-1}[0,1])\,(\forall z\in W_2^1[0,1])\qquad
	\langle Jy,z\rangle_{W_2^1[0,1]}=\langle y,z\rangle,
\]
представляет собой равномерно положительный ограниченный оператор в гильбертовом
пространстве \(W_2^1[0,1]\). Введя обозначение \(S\rightleftharpoons [JT_{\rho}(
-\varepsilon)]^{-1/2}\), находим теперь, что последовательности \(\{\lambda_n+
\varepsilon\}_{n=0}^{\infty}\) и \(\{\mu_n+\varepsilon\}_{n=0}^{\infty}\) образуют,
соответственно, спектры пучков
\begin{gather*}
	1-\lambda SJ\rho S,\\
	(1+K)-\lambda SJ\rho S,
\end{gather*}
где \(K\) "--- некоторый неотрицательный оператор конечного ранга
[\ref{par:sfkt}.\ref{pt:2.2}\,\eqref{eq:rp3}]. Согласно \cite[6.2]{VSh2},
это означает, что последовательность частичных сумм числового ряда
\[
	\sum\limits_{n=0}^{\infty}\ln\left(\dfrac{\mu_n+\varepsilon}{%
	\lambda_n+\varepsilon}\right)
\]
является неубывающей и ограниченной. Отсюда, с учётом асимптотики
\ref{par:1}.\ref{pt:1}\,\eqref{eq:mes_asymp} и неравенства \(D<1\)
[\ref{par:sfkt}.\ref{pt:sfkt}], немедленно вытекает искомое.
\end{proof}

\subsection
Перейдём теперь к доказательству собственно утверждения \ref{sing}. Через
\(\{\lambda_n\}_{n=0}^{\infty}\) при этом мы, как и при доказательстве утверждения
\ref{3.2.1}, будем обозначать последовательность занумерованных в порядке
возрастания собственных значений задачи \ref{pt:3.1}\,\eqref{eq:3.0},
\ref{pt:3.1}\,\eqref{eq:3.00}.

Заметим [\ref{par:1}.\ref{pt:1}\,\eqref{eq:mes_asymp}], что при почти всех
\(t\in [0,\nu]\) выполняется равенство \(\sigma(t)=\lim\limits_{k\to\infty}
\sigma_k(t)\), где положено \(\sigma_k(t)\rightleftharpoons\varkappa^{-k}
N(e^{k\nu+t})\). Заметим также, что независимо от выбора индекса \(k\in\mathbb N\)
для всех \(t\in [0,\nu]\), удовлетворяющих при некотором \(n\in\mathbb N\)
неравенствам
\[
	\lambda_{\varkappa (n+1)-1}<e^{(k+1)\nu+t}<\lambda_{\varkappa (n+1)},
\]
значения функций \(\sigma_k\) и \(\sigma_{k+1}\) совпадают
[\ref{par:2}.\ref{3.2.2}], а для почти всех прочих \(t\in [0,\nu]\) различаются
не более чем на \(\varkappa^{-k}\). При этом последовательность частичных сумм ряда
\[
	\sum\limits_{n=1}^{\infty}|\ln\lambda_{\varkappa (n+1)-1}-
		\ln\lambda_{\varkappa n}|
\]
является ограниченной [\ref{par:2}.\ref{3.2.2}, \ref{par:2}.\ref{3.2.3},
\ref{3.2.1}], а потому при \(k\to\infty\) мера множеств
\[
	\{t\in [0,\nu]\::\:\sigma_{k+1}(t)\neq\sigma_k(t)\}
\]
имеет асимптотику \(o(1)\). Соответственно, справедливы оценки
\begin{gather*}
	\|\sigma_{k+1}-\sigma_k\|_{L_2[0,\nu]}=o(\varkappa^{-k}),\\
	\intertext{а тогда и вытекающая из них асимптотика}
	\|\sigma_k-\sigma\|_{L_2[0,\nu]}=o(\varkappa^{-k}).
\end{gather*}
С учётом того факта, что число точек разрыва функций \(\sigma_k\) заведомо
[\ref{par:1}.\ref{pt:1}\,\eqref{eq:mes_asymp}] допускает при \(k\to\infty\)
оценку класса \(O(\varkappa^k)\), это и означает [\ref{par:5}.\ref{asymp}]
справедливость доказываемого утверждения.


\section{Примеры}\label{par:4}
\subsection
Все таблицы из настоящего параграфа содержат данные, относящиеся к уравнению
Штурма--Лиувилля, весовой функцией в котором выступает обобщённая производная
канторовой лестницы.

\begin{table}[t]
\begin{center}
\begin{tabular}{|r|r@{\;}c@{\;}l|r@{\;}c@{\;}l|r@{\;}c@{\;}l|}
\hline
{\(n\)}&\multicolumn{3}{|c|}{\(\lambda_n\)}&
\multicolumn{3}{|c|}{\(6\lambda_n\)}&\multicolumn{3}{|c|}{\(\lambda_{2n}\)}\\ \hline
1&\(7,0974\)&\(\pm\)&\(10^{-4}\)&\(42,584\)&\(\pm\)&\(10^{-3}\)&
	\(42,584\)&\(\pm\)&\(10^{-3}\)\\
2&\(42,584\)&\(\pm\)&\(10^{-3}\)&\(255,51\)&\(\pm\)&\(10^{-2}\)&
	\(255,51\)&\(\pm\)&\(10^{-2}\)\\
3&\(61,344\)&\(\pm\)&\(10^{-3}\)&\(368,06\)&\(\pm\)&\(10^{-2}\)&
	\(368,06\)&\(\pm\)&\(10^{-2}\)\\
4&\(255,51\)&\(\pm\)&\(10^{-2}\)&\(1533,0\)&\(\pm\)&\(10^{-1}\)&
	\(1533,0\)&\(\pm\)&\(10^{-1}\)\\
5&\(272,98\)&\(\pm\)&\(10^{-2}\)&\(1637,9\)&\(\pm\)&\(10^{-1}\)&
	\(1637,9\)&\(\pm\)&\(10^{-1}\)\\
6&\(368,06\)&\(\pm\)&\(10^{-2}\)&\(2208,4\)&\(\pm\)&\(10^{-1}\)&
	\(2208,4\)&\(\pm\)&\(10^{-1}\)\\
7&\(383,55\)&\(\pm\)&\(10^{-2}\)&\(2301,3\)&\(\pm\)&\(10^{-1}\)&
	\(2301,3\)&\(\pm\)&\(10^{-1}\)\\
8&\(1533,0\)&\(\pm\)&\(10^{-1}\)&\(9198,2\)&\(\pm\)&\(10^{-1}\)&
	\(9198,2\)&\(\pm\)&\(10^{-1}\)\\
9&\(1548,0\)&\(\pm\)&\(10^{-1}\)&\(9288,3\)&\(\pm\)&\(10^{-1}\)&
	\(9288,3\)&\(\pm\)&\(10^{-1}\)\\
\hline
\end{tabular}
\end{center}

\vspace{0.5cm}
\caption{Оценки первых собственных значений задачи Неймана для случая
\(\varkappa=2\), \(a=b=1/3\).}
\label{tab:1}
\end{table}
В таблице \ref{tab:1} представлены результаты численных расчётов для первых девяти
положительных собственных значений задачи Неймана. Данные таблицы иллюстрируют
утверждение \ref{par:2}.\ref{3.2.2}.

\begin{table}[t]
\begin{center}
\begin{tabular}{|r|r@{\;}c@{\;}l|r@{\;}c@{\;}l|r@{\;}c@{\;}l|}
\hline
{\(n\)}&\multicolumn{3}{|c|}{\(\mu_n\)}&
\multicolumn{3}{|c|}{\(6\mu_n\)}&\multicolumn{3}{|c|}{\(\lambda_{2n+1}\)}\\ \hline
0&\(3,2983\)&\(\pm\)&\(10^{-4}\)&\(19,790\)&\(\pm\)&\(10^{-3}\)&
	\(19,790\)&\(\pm\)&\(10^{-3}\)\\
1&\(12,558\)&\(\pm\)&\(10^{-3}\)&\(75,349\)&\(\pm\)&\(10^{-3}\)&
	\(75,349\)&\(\pm\)&\(10^{-3}\)\\
2&\(48,946\)&\(\pm\)&\(10^{-3}\)&\(293,67\)&\(\pm\)&\(10^{-2}\)&
	\(293,67\)&\(\pm\)&\(10^{-2}\)\\
3&\(65,832\)&\(\pm\)&\(10^{-3}\)&\(394,99\)&\(\pm\)&\(10^{-2}\)&
	\(394,99\)&\(\pm\)&\(10^{-2}\)\\
4&\(261,64\)&\(\pm\)&\(10^{-2}\)&\(1569,8\)&\(\pm\)&\(10^{-1}\)&
	\(1569,8\)&\(\pm\)&\(10^{-1}\)\\
\hline
\end{tabular}
\end{center}

\vspace{0.5cm}
\caption{Оценки первых собственных значений задачи \(y^{[1]}(0)-2y(0)=0\),
\(y^{[1]}(1)+2y(1)=0\) для случая \(\varkappa=2\), \(a=b=1/3\).}
\label{tab:2}
\end{table}
В таблице \ref{tab:2} представлены результаты численных расчётов для первых пяти
собственных значений граничной задачи \(y^{[1]}(0)-2y(0)=y^{[1]}(1)+2y(1)=0\).
Четвёртый столбец содержит собственные значения граничной задачи \(y^{[1]}(0)-
6y(0)=y^{[1]}(1)+6y(1)=0\). Данные таблицы иллюстрируют утверждение
\ref{par:2}.\ref{3.2.3}.

\begin{table}[t]
\begin{center}
\begin{tabular}{|r|r@{\;}c@{\;}l|r@{\;}c@{\;}l|r@{\;}c@{\;}l|}
\hline
{\(n\)}&\multicolumn{3}{|c|}{\(\mu_n\)}&
\multicolumn{3}{|c|}{\(6\mu_n\)}&\multicolumn{3}{|c|}{\(\lambda_{2n+1}\)}\\ \hline
0&\(1,1829\)&\(\pm\)&\(10^{-4}\)&\(7,0974\)&\(\pm\)&\(10^{-4}\)&
	\(7,0974\)&\(\pm\)&\(10^{-4}\)\\
1&\(10,224\)&\(\pm\)&\(10^{-3}\)&\(61,344\)&\(\pm\)&\(10^{-3}\)&
	\(61,344\)&\(\pm\)&\(10^{-3}\)\\
2&\(45,497\)&\(\pm\)&\(10^{-3}\)&\(272,98\)&\(\pm\)&\(10^{-2}\)&
	\(272,98\)&\(\pm\)&\(10^{-2}\)\\
3&\(63,925\)&\(\pm\)&\(10^{-3}\)&\(383,55\)&\(\pm\)&\(10^{-2}\)&
	\(383,55\)&\(\pm\)&\(10^{-2}\)\\
4&\(258,01\)&\(\pm\)&\(10^{-2}\)&\(1548,0\)&\(\pm\)&\(10^{-1}\)&
	\(1548,0\)&\(\pm\)&\(10^{-1}\)\\
\hline
\end{tabular}
\end{center}

\vspace{0.5cm}
\caption{Оценки первых собственных значений задачи \(y^{[1]}(0)=0\),
\(y^{[1]}(1)+2y(1)=0\) для случая \(\varkappa=2\), \(a=b=1/3\).}
\label{tab:3}
\end{table}
В таблице \ref{tab:3} представлены результаты численных расчётов для первых пяти
собственных значений граничной задачи \(y^{[1]}(0)=y^{[1]}(1)+2y(1)=0\). Четвёртый
столбец содержит взятые из таблицы \ref{tab:1} собственные значения задачи
Неймана. Данные таблицы иллюстрируют утверждение \ref{par:2}.\ref{3.2.100}.

Для получения приведённого иллюстративного материала нами была использована
вычислительная методика, описанная в работе \cite{V}.



\begin{thebibliography}{99}
\bibitem{SV} M.~Solomyak, E.~Verbitsky. \emph{On a spectral problem related to
self-similar measures}// Bull. London Math.~Soc. "--- 1995. "--- V.~27, \No~3. "---
P.~242--248.

\bibitem{VSh1} А.~А.~Владимиров, И.~А.~Шейпак. \emph{Самоподобные функции
в пространстве \(L_2[0,1]\) и задача Штурма--Лиувилля с сингулярным индефинитным
весом}// Матем. сборник. "--- 2006. "--- Т.~197, \No~11. "--- С.~13--30.

\bibitem{Naz} А.~И.~Назаров. \emph{Логарифмическая асимптотика малых уклонений
для некоторых гауссовских процессов в \(L_2\)-норме относительно самоподобной
меры}// Записки науч.~семинаров ПОМИ. "--- 2004. "--- Т.~311. "--- С.~190--213.

\bibitem{VSh2010} А.~А.~Владимиров, И.~А.~Шейпак. \emph{Асимптотика собственных
значений задачи Штурма–Лиувилля с дискретным самоподобным весом}//
Матем. заметки. "--- 2010. "--- Т.~88, \No~5. "--- С.~662--672.

\bibitem{BASh2009} А.~А.~Шкаликов, Ж.~Бен Амара. \emph{Осцилляционные теоремы
для задач Штурма-Лиувилля с потенциалами-распределениями}// Вестник МГУ. Серия~1:
Матем., мех. "--- 2009. "--- \No~3. "--- С.~40--43.

\bibitem{V2} А.~А.~Владимиров. \emph{К осцилляционной теории задачи Штурма--Лиувилля
с сингулярными коэффициентами}// Журнал выч.~матем. и матем.~физ. "--- 2009. "---
Т.~49, \No~9. "--- С.~1609--1621.

\bibitem{BNT} E.~J.~Bird, S.-M.~Ngai, A.~Teplyaev. \emph{Fractal laplacians
on the unit interval}// Ann.~Sci.~Math. Qu\'ebec. "--- 2003. "--- V.~27,
\No~2. "--- P.~135--168.

\bibitem{VSh2004} А.~А.~Владимиров, И.~А.~Шейпак. \emph{Особенности условий
Неймана в задаче Штурма--Лиувилля с сингулярным весом}// Междунар. конф. "<Дифф.
уравнения и смежные вопросы">, посвящённая 103-летию И.~Г.~Петровского. Сборник
тезисов. М.: МГУ, 2004. С.~238--239.

\bibitem{Vl} А.~А.~Владимиров. \emph{О сходимости последовательностей
обыкновенных дифференциальных операторов}// Матем.~заметки. "--- 2004. "---
Т.~75, \No~6. "--- С.~941--943.

\bibitem{Na} М.~А.~Наймарк. \emph{Линейные дифференциальные операторы}. М.: Наука,
1969.

\bibitem{Sh} И.~А.~Шейпак. \emph{О конструкции и некоторых свойствах
самоподобных функций в пространствах \(L_p[0,1]\)}// Матем. заметки. "--- 2007. "---
Т.~81, \No~6. "--- С.~924--938.

\bibitem{SaSh} А.~М.~Савчук, А.~А.~Шкаликов. \emph{Операторы Штурма--Лиувилля
с потенциалами-распределениями}// Труды~Моск.~матем. общества. "--- 2003. "---
Т.~64. "--- С.~159--212.

\bibitem{LSY} P.~Lancaster, A.~Shkalikov, Qiang~Ye. \emph{Strongly
definitizable linear pencils in Hilbert space}// Integr. Equat. Oper. Th. "---
1993. "--- V.~17. "--- P.~338--360.

\bibitem{Saks:1949} С.~Сакс. \emph{Теория интеграла}. М.: ИЛ, 1949.

\bibitem{VSh2} А.~А.~Владимиров, И.~А.~Шейпак. \emph{Асимптотика собственных
значений задачи высшего чётного порядка с дискретным самоподобным весом}//
\texttt{arXiv:1009.5335}.

\bibitem{V} А.~А.~Владимиров. \emph{О вычислении собственных значений задачи
Штурма--Лиувилля с фрактальным индефинитным весом}// Журнал выч.~матем.
и матем.~физ. "--- 2007. "--- Т.~47, \No~8. "--- С.~1350--1355.
\end{thebibliography}
\end{document}